\pdfoutput=1
\documentclass[runningheads]{llncs}

\usepackage[T1]{fontenc}
\usepackage[utf8]{inputenc}
\usepackage{amsmath}
\usepackage{newtxtext}
\usepackage[varvw]{newtxmath}
\usepackage{microtype}
\usepackage{booktabs}
\usepackage{graphicx}
\usepackage{tabularx}
\usepackage{xurl}
\usepackage{enumitem}
\usepackage[colorlinks=true,linkcolor=black,citecolor=black,urlcolor=blue]{hyperref}
\newcommand{\HH}{\mathbb H}
\newcommand{\R}{\mathbb R}
\newcommand{\Z}{\mathbb Z}

\newcommand{\cL}{\mathcal L}
\newcommand{\cA}{\mathcal A}

\newcommand{\Sym}{\operatorname{Sym}}

\newcommand{\ip}[2]{\langle #1,#2\rangle}

\makeatletter
\newcommand{\bbar}[1]{{\mathpalette\bbar@{#1}}}
\newcommand{\bbar@}[2]{%
  \setbox\z@\hbox{$\m@th#1#2$}%
  \vbox{\offinterlineskip
    \hbox to\wd\z@{\hss\vrule height 1pt depth\z@ width .75\wd\z@\kern-.05\wd\z@\hss}%
    \kern.8pt\box\z@}}
\makeatother

\title{Latin cubes of order $4$ and $4320$-point quaternionic constructions in $\R^{16}$}
\titlerunning{Latin cubes of order $4$}
\author{Rustem Takhanov}
\authorrunning{R. Takhanov}
\institute{School of Computing and Artificial Intelligence,  Nazarbayev University\\
 Nazarbayev University Research Administration\\
\email{rustem.takhanov@nu.edu.kz}}
\hypersetup{pdftitle={Latin cubes of order $4$ and $4320$-point quaternionic constructions in $\R^{16}$},pdfauthor={Rustem Takhanov}}

\begin{document}
\maketitle

\begin{abstract}
Possibly the most enduring lower-bound record for kissing numbers in dimensions greater than $8$ (excluding $24$) is
$$
\kappa_{16}\geq 4320.
$$
This lower bound is attained by a shell of the famous Barnes--Wall lattice, which was discovered in 1959.
Besides this configuration, at least one alternative $4320$-point kissing arrangement is currently known.

Recently, it was noticed that all these configurations can be described in terms of the so-called quaternionic construction,
a general type of construction that was instrumental in explaining the geometry of the $841$-point kissing arrangement in ${\mathbb R}^{12}$~\cite{T}.

We give a simple generalization of quaternionic constructions in ${\mathbb R}^{16}$ obtained from record configurations. They are parameterized by three Latin cubes of order $4$, six permutations of $3$-element set, and are all infinitesimally jammed. Computer-aided counting shows that there are $4\,072\,728\,621$ isometry classes in this family of $4320$-point kissing arrangements.
\end{abstract}
\section{Introduction}\label{sec:introduction}

The \emph{kissing number} $\kappa_n$ is the largest number of congruent balls
with pairwise disjoint interiors that can simultaneously touch another ball
of the same radius in $\R^n$. Equivalently, it is the largest cardinality of
a set $X$ of unit vectors in $\R^n$ satisfying
\[
\ip xy\leq \frac12 \qquad\text{for all distinct }x,y\in X.
\]
Indeed, with the central unit ball centered at the origin, the other centers
lie at distance $2$ from the origin, and dividing them by $2$ gives precisely
this inner-product condition. We call such a set $X$ a \emph{kissing arrangement};
see~\cite{SPLAG} for background.

Despite its elementary formulation, the kissing number problem has been
solved exactly only in dimensions $1,2,3,4,8$, and $24$. 
The first two cases are elementary. The three-dimensional case was settled
by Sch\"utte and van der Waerden~\cite{SW}, and the four-dimensional case by
Musin~\cite{Musin}. In dimensions $8$ and $24$, the minimal vectors of the
$E_8$ and Leech lattices attain the sharp upper bounds proved independently
by Levenshtein~\cite{Lev} and Odlyzko and Sloane~\cite{OS}.

Among the unresolved cases, dimensions $5,6$, and $7$ are the most natural
candidates for further exact solutions. Their best known lower bounds 
are attained by the classical root systems $D_5$, $E_6$, and $E_7$,
respectively. These constructions go back to the nineteenth-century work
of Korkine and Zolotareff~\cite{KZ}. Their persistence for more than a century
makes optimality a compelling conjecture, and a strict improvement over
any of these values would be surprising. 

Beyond dimension $8$, a similar historical argument singles out dimension
$16$. Among the unresolved cases $8<n\leq24$, the oldest surviving
lower-bound record is
\[
\kappa_{16}\geq4320,
\]
attained by the normalized minimal vectors of the Barnes--Wall lattice,
discovered in 1959~\cite{BW}. Having resisted improvement for more than six
decades, this bound makes dimension $16$ a natural next candidate for an
exact solution beyond the classical root-system cases.

The present paper shows, however, that the geometry of $4320$-point kissing
arrangements is much richer than the classical lattice example alone
suggests. Building on the quaternionic description in~\cite{T}, we
construct a family of $4320$-point arrangements in $\R^{16}$ parameterized
by Latin cubes of order $4$ and permutations matching the quaternionic
levels. We report a computer-aided enumeration yielding exactly
$4\,072\,728\,621$ isometry classes within this family. Thus, even within
this structured class of constructions, the record value is attained by
billions of geometrically distinct arrangements.

This abundance does not by itself imply that a larger kissing arrangement
exists, especially taking into account that each of those new arrangements is infinitesimally jammed. 
It does show that the longevity of the Barnes--Wall record lower bound should
not be mistaken for a complete understanding of the arrangements attaining
it. Understanding
the diversity of arrangements at the current record is a natural part of
its investigation.

\section{The Latin-cube construction}\label{sec:construction}

\subsection{Preliminaries on the quaternionic description of $24$-cells}
Let $\HH$ be the algebra of real quaternions with basis $1,i,j,k$, conjugation $q\mapsto\bar q$ and inner product
$\ip pq=\operatorname{Re}(p\bar q)$; the norm $|q|^2=q\bar q$ is multiplicative, so left and right multiplications by
unit quaternions are orthogonal maps of $\HH\cong\R^4$. We write $\HH^4\cong\R^{16}$ and denote by
$\iota_u\colon\HH\to\HH^4$ ($u=0,1,2,3$) the inclusion of the $u$-th \emph{block}. Put
\[
Q_8=\{\pm1,\pm i,\pm j,\pm k\},\qquad \omega=\tfrac{1+i+j+k}{2},\qquad \tau=\tfrac{1+i}{\sqrt2},
\]
and, for $r\in\Z/3\Z$,
\[
C_r=\omega^rQ_8,\qquad D_r=\tau\omega^rQ_8,\qquad C=C_0\sqcup C_1\sqcup C_2,\qquad D=D_0\sqcup D_1\sqcup D_2 .
\]
We use the following standard facts (see \cite[\S II]{T} or \cite[Ch.~3]{CS}). The element $\omega$ normalises $Q_8$ and
$\omega^3=-1$, so $C$ is a group of order $24$ (the binary tetrahedral group); explicitly
$C=Q_8\cup\{\tfrac12(\pm1\pm i\pm j\pm k)\}$, $C_1$ consists of the elements $\tfrac12(\pm1\pm i\pm j\pm k)$ with an even number of minus
signs and $C_2$ of those with an odd number. Moreover $\tau$ normalises $C$, $D=\tau C=C\tau$ is the set of the $24$
quaternions $(\pm e\pm e')/\sqrt2$ with $e\neq e'$ in $\{1,i,j,k\}$, and $2O=C\cup D$ is a group of order $48$ with
$CD=DC=D$, $DD=C$ (the binary octahedral group). The sets $C$ and $D$ are the vertex sets of two regular $24$-cells.

\begin{lemma}[{\cite[Lemma A.1]{T}}]\label{lem:ip} We have
\begin{enumerate}[label=(\roman*)]
\item If $x,y$ lie in the same $C_r$, or in the same $D_r$, then $\ip xy\in\{0,\pm1\}$, and $\ip xy\neq0$ iff $y=\pm x$.
\item If $x\in C_r$, $y\in C_s$ with $r\neq s$, or $x\in D_r$, $y\in D_s$ with $r\ne s$, then $\ip xy=\pm\frac12$.
\item If $x\in C$ and $y\in D$, then $\ip xy\in\{0,\pm\frac1{\sqrt2}\}$.
\end{enumerate}
\end{lemma}
\begin{proof}
(i) Each $C_r$, $D_r$ has the form $gQ_8$ with $g\in\HH$, $|g|=1$, and $Q_8$ consists of an orthonormal basis vectors and their negations.
(ii) An element of $Q_8$ and a half-integral element have inner product $\pm\frac12$; elements of $C_1$ and $C_2$ have sign
patterns of different parity, which differ in an odd number of places, so their inner product is
$\frac14(4-2\cdot\mathrm{odd})=\pm\frac12$. The statement for $D_r$, $D_s$ follows from their isometry to $C_r$, $C_s$.
(iii) Elements of $D$ have exactly two nonzero coordinates $\pm\frac1{\sqrt2}$; the dot product with an element of $Q_8$ gives
$0,\pm\frac1{\sqrt2}$ and the dot product with a half-integral element gives $\frac1{2\sqrt2}(\pm1\pm1)$.
\end{proof}

By (i)--(ii), the three sets $C_r$ are exactly the equivalence classes w.r.t. the equivalence relation ``$x=\pm y$ or $x\perp y$'' on $C$, and
likewise the $D_r$ on $D$; we call $C_0,C_1,C_2$ the \emph{$C$-levels} and $D_0,D_1,D_2$ the \emph{$D$-levels}.

Let $K=Q_8/\{+1,-1\}=\{\bar1,\bar i,\bar j,\bar k\}$, a Klein four-group, and let $q\colon K\to\{1,i,j,k\}$ be the map sending each class to its representative with positive sign, $q(\pm a)=a$; we write $q_a=q(a)$. For $r\in\Z/3\Z$ the four vectors $\omega^rq_a$ ($a\in K$) form an orthonormal basis of $\HH$ and
\[
C_r=\{\pm\omega^r q_a\mid  a\in K\}.
\]
We call $\{\pm\omega^rq_a\}$ the \emph{axis of level $r$ with label $a$}; thus $C$ has twelve axes, four in each level.

\subsection{Latin cubes of order $4$}
\begin{definition}
A \emph{Latin cube} (of order $4$, over $K$) is a set $\cA\subset K^4$ of $64$ words such that any two distinct words
differ in at least two coordinates. We write $\cL$ for the set of Latin cubes.
\end{definition}

\begin{lemma}\label{lem:latin}
For $\cA\subset K^4$ the following are equivalent:
\begin{enumerate}[label=(\alph*)]
\item $\cA$ is a Latin cube;
\item for every $u\in\{0,1,2,3\}$ the map $\cA\to K^3$ deleting the $u$-th coordinate is a bijection;
\item $\cA=\{(a,b,c,F(a,b,c))\}$ for a map $F\colon K^3\to K$ which is bijective in each argument when the other two are
fixed.
\end{enumerate}
In particular, for each coordinate $u$ and each $b\in K$ exactly $16$ words of $\cA$ satisfy $a_u=b$, and for any two
coordinates $u\neq v$ and $b,b'\in K$ exactly $4$ words satisfy $(a_u,a_v)=(b,b')$.
\end{lemma}
\begin{proof}
((a)$\Rightarrow $(b)) If distinct words cannot have three equal components, deleting a coordinate is injective, hence bijective as $|\cA|=64=|K^3|$.

((b)$\Rightarrow $(a))  Conversely, two words differing in exactly one place $u$ have the same image after deleting $u$. 

The equivalence of (b) and
(c) is straightforward.
\end{proof}
McKay and Wanless~\cite{MW}
give $64$ reduced Latin cubes of order $4$. Consequently, the number
of labeled Latin cubes of order $4$ is $|\mathcal{L}| = 4!\,(3!)^2 \cdot 64 = 55\,296$.

\subsection{Kissing arrangements based on Latin cubes}
\begin{definition}\label{def:X}
Let $\pi=(\pi_{uv})_{0\le u<v\le3}$ with $\pi_{uv}\in S_3=\Sym\{0,1,2\}$, and let $\cA=(\cA_0,\cA_1,\cA_2)\in\cL^3$.
Define $X(\pi,\cA)=E\cup\bigcup_{u<v}Y_{uv}\cup M_0\cup M_1\cup M_2\subset\HH^4$ by
\begin{align*}
E&=\{\iota_u(c)\mid \ 0\le u\le3,\ c\in C\},\\
Y_{uv}&=\Bigl\{\tfrac1{\sqrt2}\bigl(\iota_u(d)+\iota_v(d')\bigr)\mid \ r\in\Z/3\Z,\ d\in D_r,\ d'\in D_{\pi_{uv}(r)}\Bigr\},\\
M_r&=\Bigl\{\tfrac12\bigl(s_0\omega^rq_{a_0},\,s_1\omega^rq_{a_1},\,s_2\omega^rq_{a_2},\,s_3\omega^rq_{a_3}\bigr)\mid \
(a_0,\cdots,a_3)\in\cA_r,\ s_i\in\{\pm1\}\Bigr\}.
\end{align*}
\end{definition}

The points of $E$, $Y_{uv}$ and $M_r$ are supported on one, two and four blocks respectively. In the notation of
\cite[Def.~VI.1]{T} this is a quaternionic construction with four axial $24$-cells; the families $Y_{uv}$ have weight
vector $(\frac1{\sqrt2},\frac1{\sqrt2})$ on the blocks $u,v$ and the families $M_r$ have weight vector
$(\frac12,\frac12,\frac12,\frac12)$.

\begin{theorem}\label{thm:kissing}
For every $\pi\in S_3^6$ and $\cA\in\cL^3$ the set $X=X(\pi,\cA)$ consists of $4320$ distinct unit vectors and
\begin{enumerate}[label=(\alph*)]
\item $\ip xy\le\frac12$ for all distinct $x,y\in X$, i.e.\ $X$ is a kissing arrangement;
\item $X=-X$;
\item $\displaystyle\sum_{x\in X}xx^{\mathsf T}=270\,I_{16}$.
\end{enumerate}
\end{theorem}

\begin{proof}
\emph{Cardinality and norms.} $|E|=4\cdot24=96$. For each pair $u<v$ there are $3\cdot8\cdot8=192$ points in $Y_{uv}$,
so $|\bigcup Y_{uv}|=6\cdot 192=1152$. For each $r$ there are $64\cdot16=1024$ points in $M_r$: different pairs $(a,s)$ give
different points because $\{\omega^rq_b\}$ is a basis. Points in different families are distinguished by their supports
or, for $M_r$ versus $M_{r'}$, by Lemma~\ref{lem:ip}(i)--(ii). Hence $|X|=96+1152+3072=4320$, and all points are unit
vectors by construction.

The statement \emph{(b)}, i.e. antipodality of $X$, follows from $C=-C$, $D_r=-D_r$, and from the fact that the signs $s_u$ in $M_r$ are free. Let us now prove the kissing property, i.e. the statement \emph{(a)}.

Suppose that $x\neq y$ are in $X$. We go through the possible cases.

The case of $(x, y)\in E\times E$ gives: inside one block this is the regular $24$-cell $C$, where $\ip xy\le\frac12$ by Lemma~\ref{lem:ip}; points
in different blocks are orthogonal.

The case of $(x, y)\in E\times Y$, i.e. $x=\iota_u(c)$ and $y=\tfrac1{\sqrt2}\bigl(\iota_u(d)+\iota_v(d')\bigr)$, gives: $\ip{\iota_u(c)}{y}=\frac1{\sqrt2}\ip cd\le\frac1{\sqrt2}\cdot\frac1{\sqrt2}=\frac12$ by
Lemma~\ref{lem:ip}(iii).

The case of $(x, y)\in E\times M$, i.e. $x=\iota_u(c)$ and $y=\tfrac12\bigl(...,s_u\omega^rq_{a_u},...\bigr)$, gives: $\ip{\iota_u(c)}{y}=\frac12\ip{c}{s_u\omega^rq_{a_u}}\le\frac12$.

The case of $(x, y)\in Y\times Y$, i.e. $x=\tfrac1{\sqrt2}\bigl(\iota_u(d)+\iota_v(d')\bigr)$ and $y=\tfrac1{\sqrt2}\bigl(\iota_{u'}(e)+\iota_{v'}(e')\bigr)$. For two points of the same $Y_{uv}$ coming from the same $r$,
$\ip xy=\frac12(\ip de+\ip{d'}{e'})$ with both terms in $\{0,\pm1\}$ (Lemma~\ref{lem:ip}(i)) and not both equal to $1$,
so $\ip xy\le\frac12$. For different $r$, both terms equal $\pm\frac12$ by Lemma~\ref{lem:ip}(ii). For $Y_{uv}$ and a
family $Y_{u'v'}$ sharing exactly one block, $\ip xy=\frac12\ip de\le\frac12$; for disjoint pairs, $\ip xy=0$.

The case of $(x, y)\in Y\times M$, i.e. $x=\tfrac1{\sqrt2}\bigl(\iota_u(d)+\iota_v(d')\bigr)$ and $y=\tfrac12\bigl(s_i\omega^rq_{a_i}\bigr)_{0\leq i\leq 3}$, gives: $\ip xy=\frac1{2\sqrt2}\bigl(\ip{d}{s_u\omega^rq_{a_u}}+\ip{d'}{s_v\omega^rq_{a_v}}\bigr)
\le\frac1{2\sqrt2}\cdot\frac2{\sqrt2}=\frac12$ by Lemma~\ref{lem:ip}(iii).

The case of $(x, y)\in M_r\times M_{r'}$, $r\neq r'$, i.e. $x=\tfrac12\bigl(s_i\omega^rq_{a_i}\bigr)_{0\leq i\leq 3}$ and $y=\tfrac12\bigl(s_i'\omega^{r'}q_{a_i'}\bigr)_{0\leq i\leq 3}$, gives: every term $\ip{s_u\omega^rq_{a_u}}{s'_u\omega^{r'}q_{a'_u}}$ equals $\pm\frac12$, so
$\ip xy\le\frac14\cdot4\cdot\frac12=\frac12$.

Finally we have the case where the latin cube property is required.
The case of $(x, y)\in M_r\times M_r$, i.e. $x=\tfrac12\bigl(s_i\omega^rq_{a_i}\bigr)_{0\leq i\leq 3}$ and $y=\tfrac12\bigl(s_i'\omega^rq_{a_i'}\bigr)_{0\leq i\leq 3}$, gives: $\ip xy=\frac14\sum_{u:\,a_u=a'_u}s_us'_u$. If $a\ne a'$ there are at most two coordinates with
$a_u=a'_u$, so $\ip xy\le\frac12$. If $a=a'$ then $s\ne s'$ and $\ip xy\le\frac14(3-1)=\frac12$.

\emph{(c)} Let $S=\sum_{x\in X}xx^{\mathsf T}$  and the partitioning of $S$ into $4\times4$ blocks is $\left[S_{uv}\right]_{0\leq u,v\leq 3}$. For $a,b\in \HH$, the expression $ab^\top$ denotes the $4\times 4$-matrix corresponding to column-row multiplication of $a,b$ treated as four-dimensional vectors. For $u\neq v$ the contributions of
$E$ and of the families $Y_{u'v'}$ with $\{u',v'\}\ne\{u,v\}$ to $S_{uv}$ are always zero; that of $Y_{uv}$ equals
$\frac12\sum_r\bigl(\sum_{d\in D_r}d\bigr)\bigl(\sum_{d'\in D_{\pi_{uv}(r)}}d'\bigr)^{\mathsf T}=0$ because
$D_r=-D_r$, and that of $M_r$ vanishes after summing over the independent signs $s_u$. Hence $S_{uv}=0$.
For the diagonal block $S_{uu}$  we have:
\begin{itemize}
\item $E$ contributes $\sum_{c\in C}cc^{\mathsf T}=6I_4$, since $C$ is a frame;
\item each of the three families $Y_{uv}$ (or $Y_{vu}$) containing the block $u$ contributes
$8\cdot\frac12\sum_{d\in D}dd^{\mathsf T}=24I_4$, every $d\in D$ occurring with $8$ partners;
\item $M_r$ contributes $\frac14\cdot16\sum_{a\in\cA_r}(\omega^rq_{a_u})(\omega^rq_{a_u})^{\mathsf T}
=4\cdot16\sum_{b\in K}(\omega^rq_b)(\omega^rq_b)^{\mathsf T}=64I_4$ by Lemma~\ref{lem:latin}.
\end{itemize}
So $S_{uu}=(6+72+192)I_4=270I_4$.
\end{proof}

\section{Infinitesimal jamming}\label{sec:infinitesimal-jamming}

We now show that every arrangement constructed above is infinitesimally
jammed. For a kissing arrangement $Z\subset {\mathbb S}^{n-1}$, an
\emph{infinitesimal motion} is a collection of vectors
$(v_x)_{x\in Z}\subset\R^n$ satisfying
\begin{equation}\label{eq:jam-motion}
\ip{x}{v_x}=0\quad(x\in Z),\qquad
\ip{v_x}{y}+\ip{x}{v_y}\leq0
\quad\text{if }\ip{x}{y}=\frac12.
\end{equation}
The arrangement is \emph{infinitesimally jammed} if every such motion has
the form $v_x=Ax$ for a single skew-symmetric matrix $A$. An infinitesimally jammed kissing arrangement is jammed, i.e. every
continuous deformation preserving the kissing inequalities consists of
congruent configurations~\cite[Theorem~2.2]{CJKT}.

\begin{theorem}\label{thm:infinitesimally-jammed}
For every $\pi\in S_3^6$ and $\cA\in\cL^3$, the kissing arrangement
$X(\pi,\cA)$ of Definition~\ref{def:X} is infinitesimally jammed in
$\R^{16}$. More strongly, writing $Y=\bigcup_{u<v}Y_{uv}$, each
$2272$-point subconfiguration $E\cup Y\cup M_t$, $t\in\{0,1,2\}$,
is already infinitesimally jammed.
\end{theorem}

\begin{proof}
Put $V_u=\iota_u(\HH)$ and $E_u=\iota_u(C)$. For reversed indices,
use $Y_{vu}=Y_{uv}$ and $\pi_{vu}=\pi_{uv}^{-1}$.
Fix $t$ and an infinitesimal motion of $E\cup Y\cup M_t$.

\emph{Fixing the points of $E$.}
For every $u<v$, the $240$ points $E_u\cup E_v\cup Y_{uv}$ form a
normalized $E_8$ root configuration in $V_u\oplus V_v$. This follows from the fact this collection is a subset of the kissing arrangement, and is therefore a kissing arrangement itself, but in an $8$-dimensional subspace $V_u\oplus V_v$. Since it has $240$ points, it can only be isometric to the normalized $E_8$ root system. 

The $E_8$ root configuration is infinitesimally jammed
\cite[Corollary~3.4]{CJKT}. Projecting the velocities onto
$V_u\oplus V_v$ therefore gives a skew-symmetric map $A^{uv}$ on
that subspace. Its diagonal block on $V_u$ is determined by the
velocities of $E_u$, because $E_u$ spans $V_u$. Consequently this diagonal block is independent
of $v$. The off-diagonal blocks from each pair are negatives of each
other's transpose. They therefore assemble into one skew-symmetric
map $A$ on $\R^{16}$. Subtracting the infinitesimal rotation $x\mapsto Ax$,
we may assume
\begin{equation}\label{eq:jam-fixed-axes}
v_x=0\quad(x\in E),\qquad
\operatorname{pr}_{V_u\oplus V_v}v_y=0\quad(y\in Y_{uv}).
\end{equation}

\emph{Bilinearity of the remaining two-block motions.}
Write $y_{uv}(d,e)=(\iota_u(d)+\iota_v(e))/\sqrt2$, and denote the
$V_w$-component of its velocity by $v_{uv}^w(d,e)$. Fix distinct
$u,v,w$ and $d\in D_r$. The points $y_{uv}(d,e)$ and $y_{uw}(d,f)$
are in contact for every
$e\in D_{\pi_{uv}(r)}$ and $f\in D_{\pi_{uw}(r)}$. Let $y=y_{uv}(d,e)$ and $z=y_{uw}(d,f)$. 
By~\eqref{eq:jam-fixed-axes}, the
$V_u$-components of both velocities vanish, so
\begin{equation}\label{eq:jam-cross}
\ip{v_y}{z}+\ip{y}{v_z}
=\frac1{\sqrt2}\Bigl(
\ip{v_{uv}^w(d,e)}{f}+\ip{v_{uw}^v(d,f)}{e}
\Bigr)\leq0.
\end{equation}
Summing over all eight choices of $e$ and all eight choices of $f$
gives zero, since each level sums to zero. Hence all these inequalities
are equalities. Comparing the equalities for $e$ and $-e$, and using
that the $f$'s span $V_w$, gives
$v_{uv}^w(d,-e)=-v_{uv}^w(d,e)$. Contacts with $Y_{vw}$ similarly
give oddness in $d$. Each level is a basis together with its negatives,
so these values extend uniquely to a bilinear map
\begin{equation}\label{eq:jam-bilinear}
\begin{split}
B_{uv,r}^w&:V_u\times V_v\longrightarrow V_w,\\
v_{uv}^w(d,e)&=B_{uv,r}^w(d,e)
\qquad(d\in D_r,\ e\in D_{\pi_{uv}(r)}).
\end{split}
\end{equation}

\emph{Tightness of the contacts with a Latin layer.}
We use an elementary property of the two $24$-cells. For every
$c\in C_t$ and every $D$-level $D_r$, exactly two vectors of $D_r$
have inner product $1/\sqrt2$ with $c$. They have the form
\begin{equation}\label{eq:jam-companion}
d^+=\frac{c+b}{\sqrt2},\qquad
d^-=\frac{c-b}{\sqrt2},\qquad b\in C_t,\quad b\perp c.
\end{equation}
Indeed, due to symmetry, we may assume that $t=0$ and $c=1=e_0$, and the points of $D$
are all $(\pm e_a\pm e_b)/\sqrt2$, and three $D$-levels
correspond to three perfect matchings of these four axes. 
As $r$ varies, the companion axes $\{\pm b\}$ are the other three axes of $C_t$.

For $m=\frac12\sum_u\iota_u(c_u)\in M_t$, where
$c_u=s_u\omega^tq_{a_u}$, write $y\sim m$ when $y\in Y$ and
$\ip{y}{m}=1/2$. There are four such neighbors $y\in Y_{uv}$ for each $D$-level $r$. 
Their sum for fixed $u,v,r$ is
$2(\iota_u(c_u)+\iota_v(c_v))$. Consequently
\begin{equation}\label{eq:jam-stress-m}
\#\{y\in Y\mid y\sim m\}=6\cdot 3\cdot 4=72,\qquad
\sum_{y\sim m}y=36m.
\end{equation}
Conversely, for each $y\in Y$,
\begin{equation}\label{eq:jam-stress-y}
\#\{m\in M_t\mid m\sim y\}=64,\qquad
\sum_{m\sim y}m=32y.
\end{equation}
In fact, each of the two supporting blocks of $y$ permits two signed
$C_t$ axes. Each pair of axis labels occurs in four words of
$\cA_t$ by Lemma~\ref{lem:latin}, and the signs in the other two
blocks are free. This gives $2\cdot2\cdot4\cdot4=64$ neighbors.
The components outside the supporting pair cancel, while those inside
sum to the corresponding components of $32y$.

Summing the contact inequalities over all $y\sim m$ and using
\eqref{eq:jam-stress-m}--\eqref{eq:jam-stress-y} gives
\begin{align*}
&\sum_{\substack{m\in M_t,\ y\in Y\\ y\sim m}}
 \bigl(\ip{v_m}{y}+\ip{m}{v_y}\bigr)\\
&\hspace{1em}=36\sum_{m\in M_t}\ip{v_m}{m}
 +32\sum_{y\in Y}\ip{v_y}{y}=0.
\end{align*}
Every summand is nonpositive, so each is zero:
\begin{equation}\label{eq:jam-tight}
\ip{v_m}{y}+\ip{m}{v_y}=0\qquad(y\sim m).
\end{equation}

\emph{Eliminating the bilinear motions.}
Fix $u<v$ and a level $r$ of $Y_{uv}$, and let $w,z$ be the other
two blocks. For a point $m\in M_t$, choose its positive neighbors
in the relevant $D$-levels as in~\eqref{eq:jam-companion}:
\[
d^{\pm}=\frac{c_u\pm b_u}{\sqrt2},\qquad
e^{\pm}=\frac{c_v\pm b_v}{\sqrt2}.
\]
Apply~\eqref{eq:jam-tight} to the four neighbors
$y_{uv}(d^\alpha,e^\beta)$ and sum with coefficients
$\alpha\beta$, where $\alpha,\beta\in\{\pm1\}$.
The terms involving $v_m$ cancel. By
\eqref{eq:jam-fixed-axes} and~\eqref{eq:jam-bilinear}, what remains is
\begin{equation}\label{eq:jam-mixed}
\ip{c_w}{B_{uv,r}^w(b_u,b_v)}
+\ip{c_z}{B_{uv,r}^z(b_u,b_v)}=0.
\end{equation}
The point obtained from $m$ by reversing the sign of $c_z$ also
belongs to $M_t$. Adding its version of~\eqref{eq:jam-mixed} gives
\[
\ip{c_w}{B_{uv,r}^w(b_u,b_v)}=0.
\]
The Latin property makes $(a_u,a_v,a_w)$ range over all of $K^3$.
For fixed $r$, the companion-axis maps in blocks $u,v$ are permutations.
Thus this equality tests $B_{uv,r}^w$ on every pair of basis vectors
against every basis vector of $V_w$, forcing $B_{uv,r}^w=0$.
All choices were arbitrary, so~\eqref{eq:jam-fixed-axes} now yields
$v_y=0$ for every $y\in Y$.

\emph{Fixing the remaining points.}
For $m\in M_t$, equation~\eqref{eq:jam-tight} gives
$\ip{v_m}{y}=0$ for all its neighbors in $Y$. These neighbors span
$\R^{16}$. Indeed, the sums of the four neighbors on a pair and a
fixed level give $2(\iota_u(c_u)+\iota_v(c_v))$; the six pair sums
span the four directions $\iota_u(c_u)$. Differences of neighbors
that vary only in block $u$ give the companion direction $\iota_u(b_u)$.
As the level varies, these are the other three basis directions in
$V_u$. Hence $v_m=0$.

This proves infinitesimal jamming of $E\cup Y\cup M_t$.
For an infinitesimal motion of the full arrangement, its restriction
to this subconfiguration is a common rotation, which we subtract.
For any point $m$ in either remaining layer, the contact inequalities
with $Y$ are now $\ip{v_m}{y}\leq0$. Their sum is
$36\ip{v_m}{m}=0$ by~\eqref{eq:jam-stress-m}, which holds for all
three layers. Every inequality is therefore an equality, and the
same spanning argument gives $v_m=0$. Thus the full arrangement is
infinitesimally jammed as well.
\end{proof}

\begin{corollary}\label{cor:latin-cube-rigidity}
Every arrangement $X(\pi,\cA)$ is isolated, up to isometry, among
$4320$-point kissing arrangements in $\R^{16}$. In particular, a
continuous deformation starting at $X(\pi,\cA)$ and preserving all
kissing inequalities consists entirely of configurations congruent
to $X(\pi,\cA)$.
\end{corollary}
\begin{proof}
Combine Theorem~\ref{thm:infinitesimally-jammed} with the implication
from infinitesimal jamming to jamming~\cite[Theorem~2.2]{CJKT}.
\end{proof}

Our next goal is to count the number of isometry classes in the family of kissing arrangements from Definition~\ref{def:X}.
\section{Counting certificate}
\label{sec:counting-certificate}

A reproducible certificate for the count of $4\,072\,728\,621$
isometry classes is available in the
\href{https://github.com/k-nic/latin-cubes}{accompanying GitHub repository},
together with the exact counting data, verification scripts,
and a detailed explanation of the method.

\section{Use of artificial intelligence}
\label{sec:artificial-intelligence}

The idea of the quaternionic construction used in this paper originates
in our earlier work~\cite{T} on the $841$-point kissing arrangement in
$\mathbb R^{12}$~\cite{TAY,TY}. In that work, the record arrangement was interpreted as a
deformation of a  $840$-point configuration together with the
addition of one point. The underlying $840$-point configuration admits a
natural description as a subset  of $(2O)^3$ with some additional scaling of quaternionic components.

We also observed in~\cite{T} that the classical optimal kissing
configurations in $\mathbb R^4$ and $\mathbb R^8$, as well as many known
configurations attaining current or past record sizes in
$\mathbb R^{12}$, $\mathbb R^{16}$, and even $\mathbb R^{20}$, admit
quaternionic descriptions. These observations motivated the search for
a general family of quaternionic kissing arrangements in $\mathbb R^{16}$.

As starting examples, we analyzed the quaternionic descriptions of the
minimal-vector shell of the Barnes--Wall lattice and the $4320$-point
arrangement constructed and called the \emph{odd kissing configuration}
by Cohn and Li~\cite{CL}. Both examples have particularly simple quaternionic
descriptions, and their common structure suggested straightforward
generalizations. We then prompted Astra of OpenAI through
a sequence of requests for increasingly general descriptions. This
process led to the formulation in terms of Latin cubes of order $4$
presented here, which was the most natural and general description we
obtained. We carefully verified the correctness of the resulting
construction and its stated properties.

We subsequently asked Astra to count the isometry classes in this
family. It produced the counting certificate, verification code, and
accompanying description discussed in the preceding section. The
certificate and scripts are provided so that the computational part of
the enumeration can be reproduced and checked independently.
\begin{credits}
\subsubsection{\ackname} The author thanks Stanislav Yun for fruitful discussions. This research has been funded by the Science Committee of the Ministry of Science and Higher Education of the Republic of Kazakhstan (Grant No. AP27510283).

\subsubsection{\discintname}
The author has no competing interests to declare that are
relevant to the content of this article.
\end{credits}

\bibliographystyle{splncs04}
\bibliography{lit}

@article{BW,
  author    = {Barnes, E. S. and Wall, G. E.},
  title     = {Some extreme forms defined in terms of abelian groups},
  journal   = {J. Austral. Math. Soc.},
  volume    = {1},
  year      = {1959},
  pages     = {47--63}
}

@misc{TAY,
  author        = {Takhanov, Rustem and Assylbekov, Zhenisbek and Yun, Stanislav},
  title         = {Structure of kissing arrangements in {$\mathbb{R}^{12}$}
                   and a place for the {$841$}st sphere},
  year          = {2026},
  eprint        = {2606.18984},
  archivePrefix = {arXiv},
  primaryClass  = {cs.IT},
  doi           = {10.48550/arXiv.2606.18984},
  url           = {https://arxiv.org/abs/2606.18984}
}

@misc{TY,
  author        = {Takhanov, Rustem and Yun, Stanislav},
  title         = {Classification of independent sets in signed {Johnson}
                   graphs and applications to kissing arrangements},
  year          = {2026},
  eprint        = {2606.03299},
  archivePrefix = {arXiv},
  primaryClass  = {cs.IT},
  doi           = {10.48550/arXiv.2606.03299},
  url           = {https://arxiv.org/abs/2606.03299}
}

@article{CJKT,
  author  = {Cohn, Henry and Jiao, Yang and Kumar, Abhinav and Torquato, Salvatore},
  title   = {Rigidity of spherical codes},
  journal = {Geometry \& Topology},
  volume  = {15},
  number  = {4},
  year    = {2011},
  pages   = {2235--2273},
  doi     = {10.2140/gt.2011.15.2235},
  eprint  = {1102.5060},
  archivePrefix = {arXiv},
  primaryClass  = {math.MG}
}

@misc{CL,
  author = {Cohn, Henry and Li, Anqi},
  title  = {Improved kissing numbers in seventeen through twenty-one dimensions},
  year   = {2024},
  note   = {arXiv:2411.04916},
  doi    = {10.48550/arXiv.2411.04916},
  url    = {https://arxiv.org/abs/2411.04916}
}

@book{CS,
  author    = {Conway, J. H. and Smith, D. A.},
  title     = {{On Quaternions and Octonions}},
  publisher = {A K Peters},
  year      = {2003}
}

@book{SPLAG,
  author    = {Conway, J. H. and Sloane, N. J. A.},
  title     = {{Sphere Packings, Lattices and Groups}},
  edition   = {Third},
  publisher = {Springer},
  year      = {1999}
}

@article{MW,
  author    = {McKay, B. D. and Wanless, I. M.},
  title     = {A census of small {Latin} hypercubes},
  journal   = {SIAM J. Discrete Math.},
  volume    = {22},
  year      = {2008},
  pages     = {719--736}
}

@misc{T,
  author    = {Takhanov, R.},
  title     = {A quaternionic construction behind $841$-point kissing arrangement in {$\mathbb{R}^{12}$}},
  year      = {2026},
  note      = {arXiv:2609.09179}
}

@article{SW,
  author    = {Sch{\"u}tte, K. and van der Waerden, B. L.},
  title     = {{Das Problem der dreizehn Kugeln}},
  journal   = {Math. Ann.},
  volume    = {125},
  pages     = {325--334},
  year      = {1953},
  doi       = {10.1007/BF01343127}
}

@article{Musin,
  author    = {Musin, O. R.},
  title     = {The kissing number in four dimensions},
  journal   = {Ann. of Math. (2)},
  volume    = {168},
  number    = {1},
  pages     = {1--32},
  year      = {2008},
  doi       = {10.4007/annals.2008.168.1}
}

@article{Lev,
  author    = {Levenshtein, V. I.},
  title     = {On bounds for packings in {$n$}-dimensional {Euclidean} space},
  journal   = {Dokl. Akad. Nauk SSSR},
  volume    = {245},
  number    = {6},
  pages     = {1299--1303},
  year      = {1979},
  note      = {English translation in Soviet Math. Dokl. 20, 417--421 (1979)}
}

@article{OS,
  author    = {Odlyzko, A. M. and Sloane, N. J. A.},
  title     = {New bounds on the number of unit spheres that can touch a unit sphere in {$n$} dimensions},
  journal   = {J. Combin. Theory Ser. A},
  volume    = {26},
  number    = {2},
  pages     = {210--214},
  year      = {1979},
  doi       = {10.1016/0097-3165(79)90074-8}
}

@article{KZ,
  author    = {Korkine, A. and Zolotareff, G.},
  title     = {Sur les formes quadratiques},
  journal   = {Math. Ann.},
  volume    = {6},
  pages     = {366--389},
  year      = {1873},
  doi       = {10.1007/BF01442795}
}

\end{document}